\documentstyle[12pt]{article}

\newcommand{\q}{{\bf q}}

\newcommand{\J} {\bf J}
\newcommand{\cm} {c_{{\bf J}_0}}
\newcommand{\gammam}{\gamma_{{\bf J}_0}}
\newcommand{\alpham}{\alpha_{{\bf J}_0}}

\newcommand{\Sm}{S_{{\bf J}_0}}
\newcommand{\Qm}{Q_{{\bf J}_0}}
\newcommand{\Omegam}{\Omega_{{\bf J}_0}}
\newcommand{\omegam }{\omega_{{\bf J}_0}}
\newcommand{\pim}{\pi_{{\bf J}_0}}

\newcommand{\w}{{\bf w}}

\newcommand{\A }{{\bf A}}
\newcommand{\n }{{\bf n}}
\newcommand{\bomega }{{\bf \omega}}

\newcommand{\R }{I \!\! R}

\newcommand{\C }{\bf   C}
\newcommand{\I }{I \!\! I}
\newcommand{\M }{I \!\! M}
\newcommand{\U }{\bf U}

\newtheorem{theorem}{Theorem}
\newtheorem{lemma}{Lemma}
\newtheorem{definition}{Definition}

\begin{document}

\title{\sc  The Geometric Phase of the Three-Body Problem}

\author{
Richard Montgomery \\Dept. of  Mathematics\\,
Univ. of California, Santa Cruz
\\Santa Cruz, CA 95064, USA\\
email:  rmont@cats.ucsc.edu
}

\maketitle

\section{Introduction and Results}

\subsection{A Reconstruction Formula}
The three-body problem concerns understanding
the motions of three point masses
travelling in space according to Newton's laws of mechanics.
The three masses form a triangle in space and Newton's equations define a
dynamical system on the space of such triangles.  The {\sc shape}
(congruence class) of the triangle is the primary variable.
Shape variables are further divided up into an overall
scale parameter  $I$, and the
similarity class of the triangle.  The similarity classes form  a
two-sphere, denoted $S$, and called the {\sc shape sphere}.  The
appearance of this sphere is  central to our whole development.   This
shape sphere also plays a central role in Moeckel's work on the
three-body problem \cite{Moeckel}.   We view the
orientation and position of the triangle in space as
 secondary
variables.
The translational part of the motion is eliminated
by the usual trick of going to center-of-mass coordinates.
Our basic question is:
 {\sc Given that the initial and final
triangles of a three-body motion are similar, what is the
  rotation $R$ which relates the two triangles} (up to scale)?

Our main result is the answer in the form of
formula (2) below.  This is an example of  {\sc reconstruction
formula:} it reconstructs part of the original dynamics
from some reduced dynamics (essentially dynamics on $S$).
We suppose that
the planes defined by the initial and final triangles and the total
angular momentum vector ${\bf J}_0$ are known.  Let $\n_0$ and  $\n_1$ be
the normal vectors to the initial and final planes.  Let $R_0$ be the
(smallest) rotation in the ${\bf J}_0 - \n_0$ plane which takes $\n_0$
to   ${\bf J}_0$ and $R_1$   the analogous   rotation in the
${\bf J}_0 - \n_1$ taking ${\bf J}_0$ to $\n_1$.
(If $\n_i$ is coincident with ${\bf J}_0$ then its $R_i$ is the identity.)
Any $R$ which takes  $\n_0$ to $\n_1$ can be written in the form:
\begin{equation}
R=  R_1 R_{{\bf J}_0} R_0
\label{eq: reconstructiona}
\end{equation}
where $R_{{\bf J}_0}$ is some rotation about
the ${{\bf J}_0}$ axis.  Let $J_0 = \| {\bf J}_0 \|$ denote the length of
the total angular momentum vector.
Our reconstruction formula is the following
integral formula for the  rotation angle $\Delta \theta$ of
our $R_{{\bf J}_0}$:
\begin{equation}
J_0 \Delta \theta = \int_0 ^{t_1} \omegam (t) dt +
\int \int _D \Omegam
\label{eq: reconstructionb}
\end{equation}
 The first  integral in this formula is called a``dynamic phase''
in the Berry phase literature (
\cite{Shapere-Wilczek}), and the second integral is called the
``geometric phase''.  The  integrand
$\omegam$  represents $J_0$ times the instantaneous
angular velocity of the moving triangle $q(t)$ about the axis ${{\bf J}_0}$.
It is    given by
\begin{equation}
\omegam  = {{{\bf J}_0}}
\cdot \I (q(t)) ^{-1} \J_0
\label{eq: instantomega}
\end{equation}
 where $\I(q)$ is the instantaneous moment
of inertia tensor of the weighted triangle $q$.  (See (\ref{eq: inertia}),
below).     The time $t_1$ of integration is the duration of the motion.
The second  integrand,
$\Omegam$ is a  closed two-form which is independent of the potential,
 given explicitly  in equation
(\ref{eq: twoforma}) below. Its geometric definition
can be found in Theorem 4 at the end of \S 2.
This two-form lives on a ``reduced
configuration space'' denoted $\Sm$ which we will describe
next.  We urge the
reader to look at the first figure which is a picture of $\Sm$
and various of its features.

The space $\Sm$ is {\bf locally} the product of two two-spheres:
$$\Sm \cong S^2 ({1 \over 2}) \times S^2 (J ) \; \; \hbox{locally}$$
where $S^2 (R)$  denotes the sphere of radius $R$.
The first sphere is the shape sphere, which naturally has radius
$1/2$.
 {\bf Globally} $\Sm$  is the non-trivial
 two-sphere bundle over the two-sphere $S$.
(There are exactly two such bundles, one being the trivial bundle.)
We will use standard
spherical coordinates $(\phi, \theta)$ on spheres $S^2 (R)$, as well as
coordinates $(z, \theta)$ where
$$z = \cos(\phi)$$
is the normalized height of a point above the equatorial circle
$\phi = \pi/2$  so that $Rz$ is the usual height.
This induces coordinates $(z_1, \theta_1, z_2 , \theta_2)$
on  $\Sm$.  We will show that
\begin{equation}
\Omegam =
J_0 \{{1 \over 2} d(z_1 z_2) \wedge d \theta_1
+ dz_2 \wedge d \theta_2 \}.
\label{eq: twoforma}
\end{equation}

The fibering spheres of $\Sm$ are represented by the second spherical
factor in the local product representation of $\Sm$ above.  Each fiber
sphere is an  instantaneous versions of the body angular momentum sphere
occuring in  the description of the  motion of a free rigid body.  Here
``instantaneous''  refers to the instantaneous shape  of the
triangle. A point on such a fibering   sphere represents the fixed
total angular momentum vector $\J_0$ {\bf viewed from a frame}
$\{\U_1, \U_2 , \U_3 \}$ {\bf attached to the moving
triangle}.  Thus the point can be
 represented by the  vector $(\J_0 \cdot
\U_1, \J_0 \cdot \U_2, \J_0 \cdot \U_3)$ of length $J_0$.
A good choice of   frame is an
orthonormal frame
  which diagonalizes
the instantaneous inertia tensor $\I$ of the triangle.
We will always take the third frame  $\U_3$ to be the normal
$\n$ to the triangle.

The height coordinate $z_1$ on the shape sphere is
proportional to the area of the triangle.
(See the appendix for a derivation of this.)  The height coordinate $z_2$
on the second sphere is the component of the total angular
momentum normal to the triangle.
It may be helpful to have explicit formulae:
$$ z_1 = 4 \sqrt {{m_1 m_2 m_3} \over {m_1 + m_2 + m_3}}{ \Delta \over
I},$$ and
$$z_2 = {1 \over J_0} {\bf J}_0 \cdot \n.$$
In the first formula the
 $m_i$ are the   masses of the three bodies.In the second formula,  the
vector   $\n$ is a unit vector normal to the triangle spanned
by  the  position vectors $\q_1, \q_2, \q_3$ of the three bodies.
The vector $\J$ is the total angular momentum of the system.
$$I =   m_1 \|  \q _1 \| ^2 +
m_2 \|  \q _2 \| ^2 + m_3 \|   \q _3 \| ^2.$$
is its polar moment of inertia.  And
$\Delta = {1 \over 2} \n \cdot (\q_2 - \q_1 ) \times (\q_3 - \q_1)$ is the
oriented area of the triangle.

The coordinate $\theta_2$, and the local splitting
of $\Sm$ into the product of spheres depends on the choice of local frame
$\{ \U_i \}$ (choice of gauge) for the moving triangle.
In our formula for the two-form $\Omegam$ we have used the
inertial eigenframe discussed above.

A three-body motion without triple collision has
natural projections to a curve in $\Sm$ and a curve in $S$.
Our assumption that the initial and final triangles are similar means
that the projected path in  $S$ is closed.  The projected
path in $\Sm$ need not be closed, but there is a canonical way to
close it.   This is depicted in the first figure as the arcs on
the fibering spheres labelled by $R_1$ and $R_2$.
The disc $D$
over which we integrate the two-form $\Omegam$ is any disc in $\Sm$
bounding the resulting closed curve.  ($\Sm$ is simply connected.)

\subsection{The planar case}

In the planar case our question  is much simpler and   has   been
solved several  times before (\cite{Iwai}, \cite{Guichardet},
\cite{Montgomery}, {\cite{Hsiang}).
A single
angle now describes the rotation relating the two similar planar
 triangles.  This angle $\Delta \theta$ is
described by the same formula as above, which simplifies
as follows.
The integrand for the dynamic phase
is given by $\omega_{{\bf J}_0}= { 1 \over I (t)} J_0$.
The integrand  $\Omega_{{\bf J}_0}$ for the geometric phase is simply
the area form $J_0 dz_1 \wedge d \theta_1$ on the shape  sphere,
normalized so that the total area of this sphere is $ J_0 2 \pi $.

The planar  three-body problem  embeds in  the spatial problem by
taking
the initial triangle to be normal to the angular momentum vector.
For such a planar motion the fiber coordinate of $\Sm$   remains fixed at
the north pole ($z_2 = 1$) because the triangle's plane remains
perpindicular to the angular momentum vector. Hence we do not have to
deal with the two-sphere bundle $\Sm$ over $S$ in this case.

\subsection{Structure of the calculation}

Our derivation of the reconstruction formula
is quite similar to our earlier  derivation
\cite{MontgomeryA}  of a  reconstruction formula for
rigid body motion.
We
construct
  a closed loop $\gamma$ and a one-form
$\alpham$ in the three-body configuration space
  such that  when we apply  Stoke's
theorem to the integral $\int_\gamma  \alpham$  we
obtain our formula.
We construct the loop $\gamma$   by concatenating  the
 three-body  motion $q(t)$  defined by Newton's
equations with several other arcs obtained  by rotating
or scaling. The one-form $\alpham$ is the component of the ``natural
mechanical connection'' $\A$ in the direction of the total angular
momentum: $\alpham - \J_0 \cdot \A$. The connection
$\A$ first appears
explicitly in  Guichardet \cite{Guichardet}.
It can be argued that it was discovered by Smale,
who certainly has a formula for our one-form  $\alpham$
 \cite{Smale}.  It was
later used up by Iwai \cite{Iwai} and
rediscovered by Shapere and Wilczek
\cite{Shapere-Wilczek}.  I used it in   \cite{Montgomery} in studying
the Falling Cat Problem.

The form $\Omegam$ of our reconstruction formula is
essentially the push down of the two-form $d \alpham$
to the quotient space of $Q$ by the two-parameter
group of rotations about $\J_0$ and scalings.  However this
quotient is a singular space.
In order to facilitate the analysis we regularize this quotient.
This is done by introducing the space
$\tilde Q$ of oriented triangles which is a branched
cover over th standard three-body configuration space $Q$.   The idea of
this regularization is due to  Hsiang \cite{Hsiang}. The space $\Sm$
on which $\Omegam$ is defined is
the quotient of $\tilde Q$ (minus the triple collision) by this
two-parameter group of rotations and scalings.

{\sc remark} The two-form $\Omegam$ is symplectic.  It
is closely related to the minimal coupling form of
Sternberg \cite{GS}.

\subsection{Structure of paper.}

In the next section we introduce some notation and
constructions basic to our goals.
Then we state the basic theorems regarding the metric and topological
structure of the quotients $S$ and, $\Sm $ some other
intermediate quotients. We also describe more carefully  the two-form
$\Omegam$.  In \S 3 we prove these theorems.  The
proofs are based on restricting to the planar three-body problem
in which case the Hopf fibration arises naturally.  In the
final section, \S 4,  we prove our reconstruction formula in the
manner outlined above. In the appendix we derive
the fact that the shape sphere is a sphere of radius $1/2$.

\subsection{Acknowledgement}
I would like to thank Wu-Yi Hsiang for suggesting this
problem, and for helpful conversation.

\section{Constructions, Notation, Theorems}

\subsection{Basic Notation}
The three-body configuration space $Q$
consists of the set of all  triples of vectors $q = (\q_1 ,
\q_2 , \q_3)$, $\q_a \in \R^3$ whose center of
 mass is at the origin: $\Sigma m_a \q_a = 0$.
The positive numbers $m_a$ are the particle masses.
We prefer to view $Q$ as the space of weighted triangles
in space.  In any case it is  is  a six-dimensional
Euclidean vector space  with squared norm
$I(q) :=  \|  q \|^2 := m_1 \|  \q _1 \| ^2 +
m_2 \|  \q _2 \| ^2 + m_3 \|   \q _3 \| ^2 $, the same function of the
introduction (the polar moment of inertia).  The
instantaneous {\sc kinetic energy} $K$ of a  path  $q(t)$
 is defined to be
$K := {1\over 2} \| \dot q \| ^2 $.
where
$$\dot q = (\dot \q_1, \dot \q_2 , \dot \q_3)$$
denotes the time derivative of the path.  The instantaneous total
energy of a motion is
$E = K + V,$
where  $V = - \Sigma_{ i \ne j} {{m_i m_j}
\over  {\| \q_i - \q_j\| }}$ is
the usual Newtonian gravitational
potential energy.   (The  choice of
potential beyond its rotational invariance plays no role in our analysis.
Any potential which is a function of the interparticle
distances alone will work.)
A {\sc three-body motion} is a solution $q(t)$ to Newton's equations:
$m_a{d ^2 \over dt ^2}  \q_a = - \nabla _{a} V$,
 $a = 1,2, 3$.
The total energy $E$ and the total angular momentum
${\J} = \Sigma m_a \q _a \times {\bf \dot q}_a$
are constant along any three-body motion.

The moment of inertial tensor $\I (q)$
of a weighted triangle $q$ is the symmetric
non-negative $3 \times 3$ matrix  defined by
\begin{equation}
\bomega \cdot \I (q) \omega = \| \bomega \times q \|^2
\label{eq: inertia}
\end{equation}
where
$ \bomega \times q = ( \omega \times \q_1 , \bomega \times \q_2,
\bomega \times \q_3)$
denotes the infinitesimal rotation of the triangle $q$
with angular velocity $\bomega$.
 $\I$ encodes that part of the
 metric   on $Q$ in the direction of
the $G$-orbits.  It satisfies the equivariance
relations:
$\I( \lambda R q) = \lambda ^2 R \I (q) R ^T $
where
$\lambda \in \R^+ $ is a homothety, or dilation
and
$R \in G = SO(3)$ is  a rotation.
(The inertia tensor  can be expressed by the formula
$\I (q) = (1)I  - \M $
where
$\M  = \Sigma m_a \q_a \otimes \q_a$
 is the standard inertia tensor,
$(1)$ is the identity matrix, and
$I = tr(\M)$
is the  polar moment of inertia.

\subsection{Oriented Triangles}

The idea in this section is due to Hsiang \cite{Hsiang}.
We will
need to choose a normal vector $\n (t)$ to the plane of our moving
triangle $q(t)$.

\begin{definition}  An oriented triangle is a pair
$\tilde q = (q,\n)$ with $q \in Q$ and $\n \in \R^3$ a unit vector normal
to the subspace of $\R^3$ spanned by the vertices
$\q_1, \q_2, \q_3$ of the triangle $q$.
(Since the center of mass is $0$ these three position vectors
are always linearly dependent.)
The set of oriented triangles will be denoted by $\tilde Q$.
\end{definition}

If the vertices of our triangle $q$ span a plane then it has two
possible orientations  $\tilde q =  ( q, \pm \n)$
where $\pm \n$ are either of the two normals to this plane.  If the
triangle   lies in a single line  an
orientation for it is any  vector $\n$ on the unit circle in the plane
orthogonal to this line.  Such triangles are called
collinear configurations.  If $q = 0$ is the
triple collision point then $\n$ is any point on the two-sphere.

\begin{lemma} $\tilde Q$ is a smooth algebraic variety.
Away from the triple collision, the natural projection $\tilde Q \to
Q$ is a branched cover,
 branched over the collinear configurations.
The rotation group $G = SO(3)$ acts freely on $\tilde Q$
away from the triple collision point
\end{lemma}

{\sc Proof}  $\tilde Q$ is the
algebraic subvariety of $Q \times S^2$ defined by the two
equations $\q_1 \cdot \n = 0$, $\q _2 \cdot \n = 0$. The differential
of these defining functions  is   full-rank everywhere.  Apply
the implicit function theorem.
The other statements are obvious. QED

It follows from this  lemma that any smooth function or covariant
tensor on $Q$ can be lifted to $\tilde Q$.
Examples are  $V, \I$ and the
Riemannian metric.  The lift will be denoted by the same symbol
as the original.  The lifted metric fails to be positive definite along
the   branching locus. The three-body equations themselves also lift to
$\tilde Q$:

\begin{lemma}  Any three-body motion $q(t)$ which does not
 consist entirely of collinear configurations has a unique
 oriented lift $\tilde q(t) \in \tilde Q$ passing through a given
initial non-collinear oriented triangle.
\end{lemma}

The proof is obvious.

\subsection{The Reduced Configuration Space, the Shape Sphere, and other
quotient spaces}

 Let
$$G(\J_0) \subset G$$
denote the one-parameter subgroup of rotations about the
angular momentum axis $\J_0$.
The quotient space $Q/G(\J_0)$ is singular, even away from the triple
collision,  due to the presence of  extra symmetry at collinear
configurations. The introduction of the space $\tilde Q$ of oriented
triangles regularizes this quotient away from the triple collision.

\begin{definition}
The {\sc reduced configuration space}
is the quotient space
$$\Qm = \tilde Q / G(\J_0), $$
with corresponding   projection denoted by
$$\pim: \tilde Q \to \Qm.$$
The {\sc reduced motion}
corresponding to  the oriented three-body motion
$\tilde q(t)$  is the projection
$ \pim(\tilde q (t))$
of this curve to the reduced configuration space.
\end{definition}

The reduced configuration space is essentially a cone over
the space  $\Sm$ which plays a central role in our
reconstruction formula.  In order to show this
  and in order to get a good picture of
both spaces, we will  also need to understand various other
quotient spaces. Set
 $$\bar Q := \tilde Q/G = \hbox{congruence classes of
oriented triangles},$$
$$Q/G = \hbox{congruence classes of triangles}$$
The action of
$ \lambda \in \R^+$ scales each triangle  by the factor
$\lambda$ and scales distances on $Q$  by this same factor.
The polar moment (squared norm) $I$ is
a $G$ -invariant function which is  homogeneous of degree
2. Let
$$S^5 = \{ I = 1 \} \subset Q$$
denote the five-sphere in the Euclidean space $Q$ and
$$\tilde S^5 \subset \tilde Q$$
be the corresponding preimage of this sphere  under the branched cover
$\tilde Q \to Q$.  Also let
$$\tilde Q^{\ast} = \tilde Q \setminus \{ 0 \},$$
and
$$Q^{\ast} = Q \setminus \{ 0 \}.$$
(The  $0$ in $\tilde Q$ represents the two-sphere of oriented
triple collisions. )
Then we have natural identifications:
$$Q^{\ast} /\R^+ \cong S^5 $$
and
$$\tilde Q^{\ast} /\R^+ \cong \tilde S^5  \subset \tilde Q.$$
The  space
$Q^{\ast}/( G \times \R ^+)$
of similarity classes of triangles is
naturally isomorphic to
$$S_+ := S^5/G \subset Q/G$$
And the space
$\tilde Q^{\ast}/( G \times \R ^+)$
of similarity classes of oriented triangles
is isomorphic to
$$S  := \tilde S^5 /G \subset \tilde Q/G.$$
Define
$$\Sm = \tilde S^5 /G(\J_0) \subset \Qm.$$

Corresponding to these spaces we have various
projections, $Q^{\ast} \to S_+$,
$\tilde Q^{\ast} \to S$, $\tilde S^5 \to \Sm$, etcetera,  denoted by
$\pi$ or $\pim$.  Note that the fibers of
$$\Sm \to S$$
are the two-spheres:
$$ \pi_{\J_0} ^{-1} (pt.) = S^2 (J_0) =
G/ G(\J_0)$$

\begin{theorem}  $S$ is a two-sphere which we call the {\sc shape
sphere}. The projection $\tilde S^5 \to S$ is the nontrivial principle
$SO(3)$ bundle over the two-sphere.
 The projection
$\Sm \to S$ is the associated  nontrivial two-sphere bundle over
the two-sphere.   The
reduced configuration space, $\Qm$
 minus the triple collision is diffeomorphic to $\Sm \times
(0, \infty)$ where the second factor is parameterized by $I$.
\end{theorem}

{\sc explanation:}  If $G$ is a connected Lie group,
then the equivalence classes of principal $G$-bundles
over an $n$-sphere are parameterized by the homotopy group
$\pi_{n-1} (G)$.  In our case  this homotopy group
is $\pi_1 (SO(3))$ which is the two-element group.
The nontrivial $G$-bundle  over $S^2$
can be realized
as follows. Identify the two-sphere with the complex projective line
$CP^1$.  Let $\gamma \to CP^1$ denote the canonical complex line-bundle
and $ \epsilon = CP^1 \times \R$ the trivial real line bundle.
Form the rank 3 real vector bundle $E= \gamma \oplus \epsilon$.
This is an oriented vector bundle with a natural fiber-inner product.
Then the nontrivial bundle, our $\tilde S^5$, is the bundle of oriented
orthonormal frames for $E$.  And the nontrivial sphere bundle,
our $\Sm$, is the unit sphere bundle of $E$.

\subsection{Metric nature of the quotients}
 $Q$, being a Euclidean space, is  a metric space.  $G$
acts on it  by isometries, so that the quotient space $Q/G$ of congruence
classes of weighted, centered triangles, inherits a metric.  This
quotient metric --sometimes called the orbital distance metric --  is
defined by declaring that the distance between two points in the quotient
is the distance between the corresponding orbits in the original space.
The dilations $\R^+$ act on $Q$ and commute with the $G$ action so they
induce an action on the quotient as well.

The {\sc cone} over a
Riemannian manifold $(X, ds^2)$,  possibly with boundary,  is the
topological space $(X \times [0, \infty) )/ (X \times \{ 0 \})$ with
associated Riemannian metric $d \lambda ^2 + \lambda^2 ds^2$.   Here
$\lambda$ is the real parameter in $[0, \infty)$.  The quotient by
``$X \times \{ 0 \}$'' means that we crush (identify) $X \times \{ 0
\}$  to  a single point, called the cone point''. The  metric tensor
and manifold structure becomes   singular
there.

\begin{theorem}  The metric space $Q/G$  of
congruence classes of triangles
is a   cone over the space $S_+$  of
similarity classes.    The cone
point corresponds to the triple collision.
$S_+$ is isometric to the closed upper hemisphere of radius
one-half.   The equator represents collinear configurations.
The dilation parameter is
$$\lambda = \sqrt{I}.$$
\end{theorem}

Replacing $Q$ by $\tilde Q$ resolves the singularity
corresponding to the collinear configurations.  The pull-back to $\tilde
Q$ of the metric on $Q$ fails to be a metric over the
collinear configurations:  it takes
no energy to rotate a line segment about its axis .  However,
dividing by the $G$-action kills these null directions so we again
get a metric on the quotient  $\tilde Q/G$.

\begin{theorem}  The metric space  $\tilde Q/G$ of
congruence classes of oriented triangles
is a   cone over the space  $S$ of
similarity classes of oriented triangles.    The cone
point corresponds to the triple collision.
 $S$ is isometric to to the  two-sphere of radius one-half.   The
equator corresponds to the collinear configurations.
The height coordinate above the equator is
$${1 \over 2} z_1 = 2 \sqrt {{{m_1 m_2 m_3} \over {m_1 + m_2 + m_3}}}{ \Delta
\over I}.$$
where  $\Delta$ is the oriented area of the triangle:
$$\Delta = {1 \over 2} \n \cdot (q_2 - \q_1) \times (\q_3 - \q_1).$$
The map $S  \to S_+$ of the sphere to the hemisphere  which
is induced by the branched cover $\tilde Q \to Q$ corresponds to
the  quotient map obtained when we identify the hemisphere  with
the quotient space obtained by identifying points of the sphere
related by the reflection about this equator.
\end{theorem}

\subsection{Relation with the symplectic reduced space}

This section is included to connect the constructions
of the previous two sections with the symplectic reduced
space construction. The results here are not used in
arriving at our formula, but may shed some light on it.

Suppose that a compact Lie group
$G$ acts freely on a manifold $Z$,
and so on $T^* Z$.  Recall that
the symplectic reduced space at the point
$J \in Lie(G)^*$  is the
sub-quotient space
$J^{-1} (\mu) /G(\J_0)$
where $J: T^* Z \to Lie(G)^*$
is the momentum map of the action, $\mu$ is a particular
fixed element of $Lie(G)^*$,  and $G(\mu) \subset G$ is its isotropy
group ( relative to the dual of the adjoint action).
Of course in our situation we write $\mu = \J_0$.
 This
symplectic reduced space
 is  diffeomorphic
to the fiber product:
$T^*(Z/G) \times_f (Z/G(\mu))$
over the quotient $Z/G$.  This follows directly
from   \cite{thesis},
or  \cite{Weinstein}, together with the fact that
$Z/G(\mu)$ is naturally identified with the co-adjoint
orbit bundle $Z \times_G (G/G(\mu))
\subset Z \times_G Lie(G)^*$ over $Z/G$.
( See also the chapter in \cite{GS} on minimal coupling.)
In a case such as the three-body problem where
the underlying dynamics can be described by a
second order equation on $Z$, the
$T^*$-part of a reduced solution curve
in the reduced space
can be recovered
from the derivative of the projection of that
 curve to
$Z/G$.  It follows that the
entire reduced curve can be recovered from its
projection to $Z/G(\mu)$.  Thus it makes
sense to call  $Z/G(\mu)$   the
{\sc reduced configuration space} at $\mu$.

\subsection{The  Connection}

Our quotients inherit
various  tensorial objects   besides metrics.
In order to describe them we proceed generally. Suppose again that we
are given a Riemannian manifold $Z$ and a group $G$ of isometries of $Z$
acting freely. From this data we can form:
\begin{itemize}
\item{} a metric on the quotient
\item{} a connection for the principal $G$-bundle $Z \to Z/G$
\item{} a fiber inner-product on the adjoint bundle $Z
\times_G Lie(G) \to Z/G$ of Lie algebras over the quotient
\end{itemize}
The metric on the quotient (orbital distance metric) we have described.

To define the connection, we define its horizontal space.
\begin{definition}  The
{\sc horizontal space} at $z \in Z$ is the orthogonal complement
at $z$ to the group orbit through $z$.  The associated connection
form $\A : TZ \to Lie(G)$ is called
the {\sc natural connection}.
 \end{definition}
The metric on $Z/G$ is a Riemannian one.  Its metric tensor
is obtained by identifying the tangent space at $\pi(z)$
with the horizontal space at $z$.
With this definition, the projection $Z \to Z/G$
has the structure of a {\sc Riemannian submersion}.

To define the fiber inner-product on the adjoint bundle, let
$$\sigma(z): Lie(G) \to T_z Z $$
denote the infinitesimal generator of the group action:
$$\sigma(z) (\omega) = {d \over {d \epsilon}} exp(\epsilon \omega) z
|_{\epsilon = 0}.$$
 Then set
$$\| \omega \|_z ^2 = \| \sigma(z) \omega \|_Z ^2.$$
This defines the fiber inner product.
 Fix a bi-invariant inner product $\cdot$ on $Lie(G)$. Using
the inner products we can construct the transpose
$$\sigma^T (z): T_z Z \to Lie(G)$$
The map $(z , \dot z) \to \sigma^T (z) (\dot z)$
is the Noether conserved quantity, or, after we identify $Lie(G)$
and $TZ$ with their duals using the inner products,
it is the momentum map $J: T^*Z \to Lie(G)^*$.
The fiber-inner product can also be written
$$ \omega \cdot \I (z) \omega = \| \omega \|_z ^2,$$
thus defining the moment of inertia tensor $\I$.  We have
$\I (z) = \sigma^T (z) \sigma (z)$.
We now have the universal formula for the connection form
associated to this situation:
$$\A (z) = \I (z) ^{-1} \circ \sigma^T (z) .$$

In our situation $\sigma^T = \J $ is the angular momentum, viewed
as a one-form with values in $\R^3$
 (the Lie algebra of $G = SO(3)$):
$${\J} = \Sigma m_a {\q _a} \times  d \q_a.$$
$\I$ is of course our moment of inertia
tensor $\I$.
The connection form is then given by
$$\A (q) = \I (q) ^{-1} \circ \J (q).$$
{\sc The horizontal space is the space of infinitesimal deformations
with zero angular momentum.}  This gives us
a physical picture of what it means for a curve to
be horizontal, and of the length of a path in
one of the quotient spaces $S$, $Q/G$, etcetera.

All of this requires a bit of care at the collinear configurations,
since the action is not free there.  This is one of the
reasons for  introducing $\tilde Q$.
 The moment of inertia tensor of an
collinear configuration has a zero eigenvector whose eigenspace is the
line through the three masses. So one might think
that $\A$ and $\omegam$ become singular at such a configuration.
But they do not.    The kernel of $\J$ and pole of
$\I^{-1}$ are represented by the rotations of the collinear configuration
configuration about its axis and they cancel.

$\A$ has a nice physical interpretation. If $q(t)$ is a three-body
motion then $\A (q(t) ( \dot q)$ is the ``best'' choice
of assignment of an angular velocity $\bf \omega$ to the motion,
given the fact that this motion need not be a rigid motion.
If it does happen to be a rigid motion,
with infinitesimal angular velocity $\bf \omega$,
then $\A (q(t) ( \dot q) = \bf \omega$.

\begin{definition}
The form $\alpham$ is the component of $\A$ along
the fixed angular momentum vector $\J_0$:
$$\alpham =  \J_0 \cdot \A.$$
It is a one-form on
on $Q^* = Q \setminus 0 $ , or by pull-back, on
$\tilde Q ^{\ast}$. \end{definition}

{\sc remark 1}
Observe that the pull-back of $\alpham$ along a
three-body motion $q(t)$ satisfies
$$q^* \alpham = \omegam (t) dt $$
where   $\omegam (t) = \omegam (q(t))$
is the instantaneous angular velocity of a three-body curve
a defined in the introduction when we were explaining
the ``dynamic term'' of our reconstruction formula.

{\sc remark 2}
Away from the triple collision, the natural projection
 $\pim:  \tilde Q
\to \Qm$ has the structure of a principal circle bundle,
the circle being $G(\J_0)$.
Its associated connection
 one-form is $\alpham$ divided by $J_0$.

\begin{theorem}
 The form $d \alpham$ pushes down to a two-form
$\Omegam$ on $\Sm$.  This is the form
described in the introduction and   given
by the explicit formula (\ref{eq: twoforma}) above .
\end{theorem}

\section{Planar configurations and proofs}

\subsection{The Hopf fibration in planar configurations}

A {\sc planar configuration} is a triangle
lying in the plane perpindicular to the angular
momentum vector $\J_0$.  If the triangle is oriented
we will take its normal to be parallel
to $\J_0$: $\J_0 \cdot \n = J_0 >0$.
The set $Q_{planar}$ of planar configurations forms
a four-dimensional Euclidean subspace of the full
configuration space.  The action of the circle
group $G(\J_0)$ on $Q_{planar}$
is isomorphic to the action of the
circle on $\C ^2$  which takes $(\zeta_1, \zeta_2)$
to $(e^{i \theta}\zeta_1, e^{i \theta}\zeta_2)$.

The intersection of the five-sphere $\{I = 1 \}$ with
$Q_{planar}$ forms a round three-sphere,
denoted either $\tilde \Sigma \subset \tilde Q$
or $\Sigma \subset Q$.  These three-spheres are
diffeomorphic under $\beta: \tilde Q \to Q$ due to the
unique choice  of orientation.  Their quotients
by $G(\J_0)$ are isometric to the two-sphere
of radius $1 \over 2$, which is the shape
sphere $S$.  Thus:
$$G(\J_0) \to \tilde \Sigma \to \Sigma/G(\J_0) =S$$
 and
$$G(\J_0) \to \Sigma \to \Sigma/G(\J_0) =S$$
are isometric as Riemannian submersions
to the standard Hopf fibration:
$$S^1 \to S^3 (1) \to S^2 ({ 1 \over 2}).$$

\subsection{Proof of theorems 1.}

Consider  the two standard local sections of the Hopf
fibration $\tilde \Sigma \to S$. The transition function
relating these sections takes values
in $G(\J_0) \subset G$.  The local sections are also
local sections for $\tilde S^5 \to S$ and as such have
the same transition function.  Restricted to the
equator the transition function represents the nontrivial
generator of the fundamental group of $G = SO(3)$
and hence $\tilde S^5 \to S$ is the nontrivial bundle.

To prove the facts regarding $\Sm$ observe that
$\tilde S^5/G(\J_0)$ is isomorphic to the
associated bundle $\tilde S^5 \times_G (G/G(\J_0))$.

{\sc proofs of theorems 2 and 3.}

Any triangle can be made to lie in the xy plane  by
a rotation so that $Q_{planar}$ is a slice for the $G$ action on
$Q$.   An oriented triangle can be made planar in a unique
way, up to rotation.
An unoriented  triangle  can be made planar in two
rotationally inequivalent  ways, the two ways being related by
reflection. In other words:  $Q/G = Q_{planar}/ O(2)$, whereas
$\tilde Q / G = Q_{planar}/ SO(2)$.  This accounts for the difference in
the two quotients. These two identifications are isometries, since
$Q_{planar}$ is totally geodesic.  The last space is $C^2/S^1$  which is
isometric to the cone over the  sphere $S^2 ({1 \over 2})$. The quotient
group $O(2)/SO(2)$ is the two-element group and  accounts for the
branched cover
 $S \to S_+$.   The derivation of
the formula for the normalized height $z_1$ can be found
in Hsiang \cite{Hsiang} and in our appendix.

 The action by homotheties commutes with
rotations so it descends to the quotient where it remains
a dilation:  $d (\lambda a , \lambda b) = \lambda d (a,
b)$.  (Here $a, b$ represent similarity classes
and $d$ is the distance function.)  Since $I$ is homogeneous
of degree 2, and since $S^5$ is defined by $I =1$, the dilation parameter
$\lambda$ equals $\sqrt I$.

{\sc proof of theorem 4.}
It follows from  the discussion of the previous section,
the above proofs
and the fact that any curve of planar triangles  has
angular momentum in the $\hat e_3$ direction,
that the restriction of $\A$ to
$\Sigma \subset Q$ is $\Gamma \hat e_3$ where
$\Gamma$ is the canonical connection for the Hopf
fibration. One can choose a local section
$s: U \subset S \to \Sigma$
for the Hopf fibration such that
$s^* \Gamma = {1 \over 2} z_1 d \theta_1$.
(The domain of this section is  the sphere
minus a geodesic arc connecting the
north and south pole.)  It follows
that
\begin{equation}
s^* \A = ({1 \over 2} z_1 d \theta_1) \hat e_3.
\label{eq: sA}
\end{equation}

Now $\A$ is ``basic'' with respect to the action of
the group $\R^+$ of dilations.  This
means that
\begin{equation}
\sigma_{\lambda}^* \A = \A
\label{eq: invariancea}
\end{equation}
\begin{equation}
{\partial \over \partial \lambda} \rfloor \A = 0
\label{eq: invarianceb}
\end{equation}
where $\sigma_{\lambda}: \tilde Q \to \tilde Q$
is homothety by $\lambda \in \R^+$
and ${\partial \over \partial \lambda} \rfloor$
denotes inner product with the infinitesimal
generator ${\partial \over \partial \lambda}$ of
homotheties.  (To see (\ref{eq: invariancea}) observe
that $\I (\lambda q) = \lambda^2 \I (q)$
and $\sigma_{\lambda}^* J =
\Sigma m_a \lambda \q_a \times d(\lambda \q_1) =
\lambda^2 \J$, and use the definition  $\A = \I^{-1}
\J$.  To see (\ref{eq: invarianceb})
observe that the angular momentum of a pure dilational
motion is $0$ which means that
${\partial \over \partial \lambda} \rfloor \A = 0$.)
Extend the section $s$ by making it constant under homothety.
Then, by the homothety invariances of $\A$,
formula (\ref{eq: sA}) still holds for this extended section.
Let
$\U$ denote the local frame induced by $s$.  It is a map to  $G = SO(3)$
 defined by writing $\tilde q = \U(\tilde q) s (\pi
(\tilde q))$.  The induced local trivialization of our
principal bundle $\tilde Q ^{\ast} \to \tilde Q ^{\ast}/G$
is then $\tilde q \mapsto (\pi (\tilde q), \U (\tilde
q))$. Note that
$\U(\tilde q)_3 = {\U}(\tilde q) (e_3) = \n$,
is
the normal vector of the oriented triangle
$\tilde q$.
Using this fact,  the transformation formula
for connections, and the fact that under
our identification of the Lie algebra of $G$ with
  $\R^3$ the
adjoint action of $G$ becomes its usual action on
$\R^3$, we see that with respect to our local
trivialization we have:
$$\A = ({1 \over 2} z_1 d \theta_1
) \n + (d U)U^{-1}$$ where
$(d U)U^{-1} = \Theta$
denotes the pull-back of the Maurer-Cartan form on $G$ by
the map $U$.

We now have
$$\alpham = {1 \over 2} (z_1 d \theta_1) J_0 z_2 + \J_0
\cdot \Theta,$$
since $\J_0 \cdot \n = J_0 z_2$.
It is well-known (see \cite{MontgomeryA}) that the two-form
$d (\J_0 \cdot \Theta)$ pushes down to
$S^2 = G/ G(\J_0)$ and that this push-down
is given by $J_0$ times the solid angle
form, $dz_2 \wedge d \theta_2$.
Thus
$d \alpham = J_0 \{ {1 \over 2} d(z_1 z_2) \wedge d \theta_1
+ dz_2 \wedge d \theta_2 \} $,
which is the claimed formula for
$\Omegam$.

{\sc Remark}.
We  think of the local frame $\U = (\U_1 , \U_2 , \U_3)$
as a moving frame attached to our triangle, chosen
so that   $\U_3 = \n$ is the triangle's normal.
Hsiang \cite{Hsiang} has shown that  the frame   which diagonalizes the
moment of inertia tensor corresponds to a
local section $s$ which also satisfies
$s^* \Gamma = {1 \over 2} z_1 d \theta_1 $. Hence  his
local trivialization is the same as the one we are
using, up to a constant rotation about the $\n$-axis.  Note that the
eigenframe  is not defined at the  north and south poles.  These
poles  correspond to the weighted triangles for
which $\I$ has double eigenvalues.  (If the masses are all
equal these are the equilateral triangles.)
A branch cut from the north to south pole is necessary,
for if we traverse a small loop encircling the fiber over
one of the poles  then we find that $(\U_1 , \U_2 , \U_3 )
\mapsto (-\U_1 , -\U_2 , \U_3 )$.

\section{Derivation of the reconstruction formula}

In this section we derive our reconstruction formula,
(\ref{eq: reconstructionb}).

\subsection{Closing the loop: a loop and a disc in $\Sm$  }

Let
$$s(t) = \pi( \tilde q (t)) \in S$$
be the curve of similarity classes represented by
our three body motion.
It is  a {\sc closed} curve on the base
two-sphere.
Let
$$c_{J}(t) = \pim ( \tilde q (t)) \in \Sm$$
denote the projection of $\tilde q (t)$ to $\Sm$.
Although $s(t)$ is closed, the reduced curve $\cm (t)$ need not be.
There is a canonical way to close it.  To see this, observe that
$\Sm \to S$ has two canonical sections.  One consists of
equivalence classes of triangles whose normals are pointing
along the $\J_0$ axis, and the other consists of those whose normals are
antiparallel to the $\J_0$ axis.  We will  these sections, or their
values at a particular similarity class,  the
``north'' and ``south'' poles.   Since $s(t)$ is closed, both  endpoints
of the curve $\cm (t)$ lie  on the same spherical fiber over
$s(0) =s(t_1)$.   On this fixed fiber draw the
two geodesic arcs from the north pole to $\cm (0)$ and from $\cm (t_1)$
back, then sandwich the reduced curve $\cm (t)$ in between.  The resulting
closed curve will be denoted $\gammam$. It is the projection of a closed
curve $\gamma(t)$ in  $\tilde Q$.  See the figure.

Finally,  $\Sm$ is simply connected so that $\gammam$ bounds some disc .
$D \subset \Sm$.

\subsection{The loop in $Q$.}

We now construct the loop $\gamma$ in $\tilde Q$
over which we integrate.   Its projection to
$\Sm$ is the loop just described above.  The loop $\gamma$ is obtained by
concatenating the dynamic curve $\tilde q(t)$ with several
group orbits, denoted $c_i (t)$ or $h(t)$.  The act of concatenating
two curves, one ending where the other begins,
is defined in the obvious manner, and will be denoted by ``$*$'' below.

To construct the group curves we will  use the following
exponential  notation for rotations. If  ${\bf v} \in \R^3$
then $exp({\bf v}) $ ,
will mean the counter-clockwise rotation about the axis spanned by
${\bf v}$ by $\| {\bf v} \|$ radians.  This is the standard Lie
theoretic exponential map if we use  the standard
identification of $\R ^3$ with the Lie algebra of the rotation
group.  If ${\bf v}$ is a unit vector, and $\theta \in \R$ then
$exp ( \theta {\bf v})$ is a rotation by $\theta$ radians
about the ${\bf v}$ axis.  Let $\n_0$
and $\n_1$ be the initial and final normal vectors
to our curve of oriented triangles, as in
assumption (2) , above.
Form unit vectors
$$\xi_0 = { 1 \over { \| \J_0 \times \n_0 \|}}\n_0 \times
\J_0$$ and
$$\xi_1 = { 1 \over { \| \n_1 \times \J_0 \|}}\J_0 \times
\n_1$$ and corresponding one-parameter
subgroups
$exp(s \xi_0)$, $exp(s \xi_1)$ of rotations.
Then:
$$R_0 = exp( \phi_2 (0) \xi_0)$$
and
$$R_1 = exp( \phi_2 (t_1) \xi_1)$$
where
$R_0, R_1$ are the rotation matrices of our reconstruction
formula (\ref{eq: reconstructiona}),$\phi_2$ is the angle in our
parameterization of $\Sm$, and $\phi_2 (t)$ is its value along
the reduced curve $\gammam (t)$:
$J \cos (\phi_2 (t)) = \J_0 \cdot \n (t)$.

Let $\tilde q(t)$ be the oriented three-body motion.
Consider the concatenation
$$\gamma = c_0 * \tilde q * c_1 * \cm *h$$
of the following curves:
$$c_0 (t) = exp(-t \xi_0) R_0 \tilde q (0), 0 \le t \le \phi_2 (0),$$
$$\tilde q (t) , 0 \le t \le t_1,$$
$$c_1 (t) = exp(t \xi_1) \tilde q (t_1),  0 \le t \le
\phi_2(t_1),$$
$$c_{\J_0} (t) = exp (-t {1 \over J_0}{\J_0}) R_1 \tilde q(t_1),
0 \le t \le \Delta \theta , $$
and
$$h(t) = e^{at} c_J (\Delta \theta).$$
The interval of definitions of the curves are
chosen so that the endpoint of one curve is the
initial point of the next and so the concatenations
are well-defined.   The constant $a$ and the time of stopping
for the final purely dilational curve $h(t)$ are chosen
so that its endpoint is the beginnning point,
$R_0 \tilde q (0)$, for $\gamma(t)$.

\subsection{Line integrals}

We have:

$$\int_\gamma \alpham =
\int_{c_1} \alpham +
\int_{\tilde q} \alpham
+ \int_{c_2} \alpham
+ \int_{c_{J}} \alpham
+ \int_{h} \alpham $$

CLAIM:
$$ \int_{c_0} \alpham = 0$$
$$ \int_{c_1} \alpham = 0$$
$$ \int_{\tilde q} \alpham = \int_0 ^{t_1} \omegam dt$$
$$\int_{\cm} \alpham = -J \Delta \theta$$
$$ \int_{h} \alpham = 0$$
The first two integrals vanish because
$$c ^* \alpham = \J_0 \cdot c^*\A
 = \J_0 \cdot \omega ds$$ whenever   $c(s) = exp(s \omega) c(0)$ is the
orbit generated by a   one-parameter subgroup of $G$.
(This follows immediately from one of  the defining properties of
connections.) Now use the fact that the infinitesimal generators $\omega =
\xi_0, \xi_1$ for the curves $c_0, c_1$ are perpindicular to $\J_0$.
To evaluate  the third integral, observe that $\cm$ is also
the orbit of a one-parameter subgroup, but its generator is the unit
vector along $\J_0$.
The vanishing of the integral over the homothety path $h$
follows immediately from the homothety invariance
 of the connection already discussed. The integrand for the dynamic path
$\tilde q (t)$ was already discussed. (See remark 1 near the end of \S
2.6.) There we noted that $q^* \alpham = \omegam dt$.

An application of Stokes' theorem and the formulae relating $d \alpham$
to $\Omegam$ now prove our reconstruction formula,
(\ref{eq: reconstructionb}).

\section{ Appendix: explicit identification of
the shape sphere}

Following the discusson of \S 3.1 it suffices to
understand the geometry of the space of similarity classes of weighted
triangles for the  planar three-body problem
We identify the plane in which the bodies move with
the complex plane.  Then we  replace the
spatial configuration space $Q$ above by the planar
configuration space $Q_{planar}$ of triples
$q = (q_1, q_2,
q_3)$ of complex numbers, subject to the constraint
$\Sigma_a m_a q_a = 0$.  The space of similarity classes
of triangles in the plane forms a two-sphere, and a
three-body motion describes a curve  $w(t)$ on this
sphere.

Let us describe the sphere $S$ of
similarity classes  explicitly.
First, we diagonalize the mass matrix
(kinetic energy) by introducting  Jacobi coordinates
$$\xi_1 = q_1 - q_3$$
and
$$\xi_2 = {- m_1 \over {m_1 + m_3}} q_1 + q_2
{- m_3 \over {m_1 + m_3}} q_3,$$
and normalized Jacobi coordinates
$$\zeta_1 = \sqrt{\mu_1} \xi_1,$$
and
$$\zeta_1 = \sqrt{\mu_2} \xi_2,$$
where the reduced masses
$\mu_i$ are defined by ${1 \over \mu_1} = {1 \over m_1} +
{1 \over m_3}$ and ${1 \over \mu_2} = {1 \over {m_1 +
m_3} } + {1 \over m_2}$.
Then
$$\Sigma m_a \| \dot q_a \| ^2 =  \| \dot \zeta_1 \|^2 +
\| \dot \zeta_2 \|^2.$$
Rotations by an angle
$\theta$ induce the transformation $(\zeta_1, \zeta_2)
\mapsto (e^{i \theta} \zeta_1, e^{i \theta} \zeta_1)$
It follows that the vector
$$\w = (w_1, w_2, w_3)$$
defined by
$$w_1 = {1 \over 2} (\| \zeta_1 \|^2 - \| \zeta_2 \|^2 ),$$
and
$$w_2 + i w_3 = \zeta_1 \bar \zeta_2$$
is invariant under rotations.  The
sphere of radius ${1 \over 2}$
$$S = \{ \w : w_1 ^2 + w_2 ^2 + w_3 ^3  = {1 \over 4} \}$$
is naturally identified with the space of similarity
classes of planar triangles.
We calculate that the height coordinate $w_3$
is
$$w_3 = {2 \over I} \sqrt{{{m_1 m_2 m_3} \over {m_1 + m_2 + m_3}}}
\Delta$$ where $\Delta$ is the (oriented) area of the
triangle $q$.
In other words,the height on $S$ represents
the   triangle's   area, $\Delta$. The normalized
height used in the body of the text is related
to this coordinate by $w_3 = {1 \over 2} z_1$.


\section{Bibliography}
\begin{thebibliography}{30}

\bibitem{Guichardet} A. Guichardet {\em On rotation and
vibration motions of molecules}
Ann. Inst. H. Poincare, Phys. Theor., {\bf 40} no.3
329-342 (1984).

\bibitem{GS} V. Guillemin and S. Sternberg, {\em Symplectic
Techniques  in Physics,} Cambridge University Press,
Cambridge, 1984.

\bibitem{Hsiang} W-Y Hsiang, {\em Geometric Study
of the Three-Body Problem, I}, UC Berkeley preprint,
PAM-620, 1994.



\bibitem{Iwai} T. Iwai, {\em
 A geometric
setting for classical molecular dynamics},Ann. Inst.
Henri Poincare, vol. 47, no. 2, pp 199-219.
(1987).



\bibitem{Levi}
M. Levi,
{\em Geometric Phases in the Motion of Rigid Bodies } ,
Archive for Rational Mechanics and Analysis,
v. 122, no. 3, 213-219, 1993.

\bibitem{MMR} J. Marsden, R. Montgomery, and T.
Ratiu,  {\em Reduction, Symmetry, and Berry's Phase in Mechanics},
Memoirs of the A.M.S.

\bibitem{MW} J. Marsden and A. Weinstein, Reduction of
Symplectic Manifolds with Symmetry, {\em Rep. Math. Phys.}
{\bf 5} (1974), 121--130.



\bibitem{Montgomery}	R. Montgomery , the appendix to:
{\em
Optimal Control of Deformable Bodies and Its Relation
to Gauge Theory} in {\bf The Geometry of Hamiltonian
Systems}, proceedings of a Workshop held in June 1988 at
M.S.R.I., T. Ratiu ed., Springer-Verlag, (1991).

\bibitem{MontgomeryA} R. Montgomery,
{\em  How Much Does a Rigid Body Rotate?},
 Am. J. Physics,v.  59, no. 5, 394-398, (1991).

\bibitem{thesis} R. Montgomery,
{\em  Bundle Theory in Mechanics},
PhD thesis, UC Berkeley, 1986.



\bibitem{Moeckel} R. Moeckel, {\em
Some qualitative features of the three-body
problem} in {\bf Hamiltonian Dynamical Systems,
proceedings, AMS Contemporary Math, Vol. 81}
, ed. K. Meyer and
C. Saari, p. 1-22, American Math. Society, (1988).

\bibitem{Shapere-Wilczek}  A. Shapere and F. Wilczek,
 {\em Geometric Phases in Physics}, World Sci.,
Singapore, (1989).

\bibitem{Smale}  S. Smale, {\em Topology and Mechanics,
I; II}, Inventiones,.v. 10, 305-311 ; v. 11, 45-64 (1970).

\bibitem{Weinstein} A. Weinstein,
{\em  A universal phase
space for a particle in a Yang-Mills field}, Lett. Math. Phys., {\bf 2},
417-420, (1978).

\end{thebibliography}
\end{document}